\documentclass[a4paper,pdftex,reqno,10pt]{amsart}

\allowdisplaybreaks[1]

\usepackage{geometry}
\usepackage{graphicx}
\usepackage{enumerate}
\usepackage{courier}
\usepackage{times}
\usepackage{amsmath,amssymb}
\usepackage{xcolor}

\usepackage[colorlinks = true, citecolor = {blue}]{hyperref}
\usepackage{cleveref}

\usepackage{lineno}
\theoremstyle{plain}
\newtheorem{Theorem}{Theorem}
\newtheorem{Proposition}[Theorem]{Proposition}
\newtheorem{Conjecture}[Theorem]{Conjecture}
\newtheorem{Corollary}[Theorem]{Corollary}
\newtheorem{Lemma}[Theorem]{Lemma}
\newtheorem{Construction}[Theorem]{Construction}

\newtheorem{Observation}[Theorem]{Observation}

\newtheorem{Question}[Theorem]{Question}
\newenvironment{Proof}
{\begin{trivlist}\item[]{{\sc Proof.}}}{\hfill{$\square$}\noindent\end{trivlist}}

\theoremstyle{definition}
\newtheorem{Definition}[Theorem]{Definition}

\theoremstyle{remark}
\newtheorem{Remark}[Theorem]{Remark}

\newcommand{\todo}[1]{{\noindent\textcolor{red}{#1}}}

\newcommand{\N}{\mathbb{N}}
\newcommand{\F}{\ensuremath{\mathbb{F}}}
\newcommand{\Fq}{\ensuremath{\F_q}}

\newcommand{\AG}{\operatorname{AG}}
\newcommand{\PG}{\operatorname{PG}}
\newcommand{\aspace}{\PG(v-1,q)}

\newcommand{\cC}{\mathcal{C}}
\newcommand{\cF}{\mathcal{F}}
\newcommand{\cK}{\mathcal{K}}
\newcommand{\cM}{\mathcal{M}}
\newcommand{\cS}{\mathcal{S}}

\begin{document}

\title{A generalization of the cylinder conjecture for divisible codes}

\author{Sascha Kurz, Sam Mattheus}
\address{Sascha Kurz, University of Bayreuth, 95440 Bayreuth, Germany \newline Sam Mattheus, Universiteit Brussel, 1050 Elsene, Belgium}
\email{sascha.kurz@uni-bayreuth.de; sam.mattheus@vub.ac.be}

\abstract{We extend the original cylinder conjecture on point sets in affine three-dimensional space to the more general framework of divisible linear codes over $\Fq$ and their classification. Through a mix of linear programming, combinatorial techniques and computer enumeration, we investigate the structural properties of these codes. In this way, we can prove a reduction theorem for a generalization of the cylinder conjecture, show some instances where it does not hold and prove its validity for small values of $q$. In particular, we correct a flawed proof for the original cylinder conjecture for $q = 5$ and present the first proof for $q = 7$. 
}}

\maketitle

\todo{
%- change notation $(q,r,v)$ to $(v,r,q)$ \\
%- add (r+1)-space implies cylinder \\
%- add strenghtening on non-existence of high weight hyperplanes \\
%- improve result for q = 3 \\
%- add counterexample based on field reduction \\
%- add ward's bound to show that $v$ cannot be unbounded \\
%- write something about classification of divisible codes (motivation from coding-theoretical point of view) \\
%$  $- rewrite introduction a little (don't start with 'Let...') \\
%- elaborate abstract \\  
%- check consistency of notation (capital letters for points, $\cS$ for divisible set of points and $\cK$ for set of points in plane in last section) \\
%- does abstract still reflect new content? \\
%- structure of the paper \\
%** start of with linear programming, obtain information on $a_0$ and $a_{q^r}$, the two most important ones. the latter as we can use induction on it \\
%** then transfer results from paper with Jan \\
%** then equivalence for (v,r,q) and (v-r+1,1,q) \\
%** then results for small q \\
%** counterexample based on field reduction \\
%** r = 1 in detail \\
%- introduce notation for $\cS(H)$ as a function counting the number of points in $H$, similarly notion of 0- and 1-point, $k$-line etc. \\
%- redo bibliography with mathscinet entries \\
%- check that all mentions of 'arc' and the usage of $\#$ are gone \\ 
%- check that all mentions of cylinder conjecture correctly use generalized
}
\section{Introduction}\label{sec_introduction}

The cylinder conjecture, as originally formulated in \cite{ball2008graph}, was motivated by direction problems in finite geometry. To be more precise, the strong version of the conjecture is formulated as follows.

\begin{Conjecture}
	Let $\cS$ be a set of $p^2$ points in $\AG(3,p)$, $p$ prime. If $\cS$ has the property that every plane is incident with $0$ modulo $p$ points of $\cS$ then $\cS$ is a cylinder, i.e. the union of $p$ parallel lines.
\end{Conjecture}

The idea to classify all sets of $p^2$ points in $\AG(3,p)$ determining few directions, is a continuation of similar results in $\AG(2,p)$, which started with the work of R\'edei and Megyesi \cite{redei1970} and Lov\'asz and Schrijver \cite{lovasz1983}. The cylinder conjecture as stated above, is of interest from a coding theoretical perspective as well. In particular, one could view this conjecture as part of a more general search for the classification of divisible codes.

A linear $[n,k]_q$ code $\cC$ is $\Delta$-divisible if all its weights are multiples of a fixed integer $\Delta$. When $\Delta$ and $q$ are coprime, the classification of $\Delta$-divisible codes over $\Fq$ is almost trivial, see for example \cite{ward2001divisible} for a general survey on $\Delta$-divisible codes, which contains this result. In this paper, we will focus on the case when $\Delta$ is a power of the field size $q$, hence looking at $q^r$-divisible codes, where $r$ is a non-negative integer. For a survey on these codes we refer to \cite{partial_spreads_and_vector_space_partitions}. In the case of the original cylinder conjecture, a positive answer would give a classification of all $[p^2,4]_p$ codes that are $p$-divisible as we will explain later on. With the classification of $q^r$-divisible codes over $\Fq$ in mind it thus makes sense to generalize the cylinder conjecture to higher dimensions and non-prime field characteristic. In order to state and motivate it, we require some terminology and a proper notion of cylinders in higher dimensions. We will therefore defer it to the next section. 

\section{Preliminaries}\label{sec_preliminaries}

The notation we will use throughout is the following. Let $V\cong{\F}_q^v$ be a $v$-dimensional vector space over the finite field {\Fq} with $q$ elements and $\aspace$ the projective space associated to it. By a $k$-space of $\aspace$ we mean a $k$-dimensional linear subspace of $V$, also using the terms points, lines, planes, and hyperplanes for $1$-, $2$-, $3$-spaces, and $(v-1)$-spaces respectively. A multiset of points in $\aspace$ will be denoted as $\cM$, while $\cS$ refers to a set of points. 

To each multiset $\cM$ of $n$ points in $\aspace$ we can assign a $q$-ary linear code $C(\cM)$ defined by its generator matrix whose $n$ columns consist of representatives of the $n$ points of $\cM$. The code $C(\cM)$ is projective if and only if $\cM$ is a set.  The weight of a codeword is the number of non-zero coordinates and as mentioned before, a code is called $q^r$-divisible if the weight of each codeword is divisible by $q^r$. 

\begin{Definition}\label{def_divisible}
	The multiset $\cM$ of $n$ points in $\aspace$ is \textbf{$q^r$-divisible} if and only if $C(\cM)$ is. Equivalently, $\cM$ is $q^r$-divisible if $|\cM \cap H| \equiv |\cM| \pmod{q^r}$ for every hyperplane $H$ of $\aspace$.
\end{Definition}

For example, the set of points of a $k$-space is $q^{k-1}$-divisible and so is the multiset of the points of a collection of $k$-spaces. The converse gives an interesting question: For which integers $l$ is each $q^{k-1}$ divisible set of $\frac{q^{k}-1}{q-1}\cdot l$ points the union of $l$ (disjoint) $k$-spaces? We can always choose $l\ge 1$ and the maximum value for $l$ implies extendability results of spreads: each set of $q^k+1-l$ disjoint $k$-spaces in $\F_q^{2k}$ can be extended to a $k$-spread. These results are mostly formulated in the language of minihypers, see e.g.\ \cite{govaerts2003particular}.

An interesting property of $q^r$-divisible sets of points is that the divisibility is preserved (to some extent) upon intersecting with subspaces, allowing for inductive arguments.

\begin{Lemma}(\cite[Lemma 7]{partial_spreads_and_vector_space_partitions})
	\label{lemma_heritable}
	Suppose that $\cM$ is a $q^r$-divisible multiset of $m$ points in $\aspace$
	and $X$ a ($v-j$)-space of $\aspace$ with $1\leq j<r$. Then the
	restriction $\cM\cap X$ is $q^{r-j}$-divisible.
\end{Lemma}
\begin{Proof}
	By induction, it suffices to consider the case $j=1$, i.e. $X=H$ is a hyperplane in $\aspace$.
	
	The hyperplanes of $H$ are the ($v-2$)-subspaces of $\aspace$ contained in $H$. Hence the assertion is equivalent to
	$|\cM\cap U|\equiv|\cM|\pmod{q^{r-1}}$ for every ($v-2$)-subspace $U\subset\aspace$. By assumption we have
	$|\cM\cap H_i|\equiv m\pmod{q^r}$ for the $q+1$ hyperplanes $H_1,\dots,H_{q+1}$ lying above $U$. This gives
	\begin{equation*}
	(q+1)m\equiv\sum_{i=1}^{q+1}|\cM\cap H_i| = q\cdot|\cM\cap U|+|\cM|\equiv q\cdot|\cM\cap U|+m\pmod{q^r}
	\end{equation*}
	and hence $m\equiv|\cM\cap U|\pmod{q^{r-1}}$, as claimed.
\end{Proof}

We can now generalize the concept of a cylinder to higher dimensional vector spaces as follows.

\begin{Definition}
	\label{def_cylinder}
	Let $r$ be a non-negative integer. An \textbf{$(r+1)$-cylinder} is a multiset of $q^{r+1}$ points in $\aspace$ that arises as the union of the points of $q$ affine $(r+1)$-subspaces $L_1\setminus F$, $\dots$, 
	$L_q\setminus F$, where the $L_i$ are $(r+1)$-spaces and $F$ is a $r$-space that is contained in all $L_i$.
\end{Definition}

We remark that our definition of a $2$-cylinder matches the definition of a cylinder in \cite{de2019cylinder} and the one stated above. By convention a $1$-cylinder is just a multiset of $q$ points. As the affine subspaces mentioned in the definition above will appear often, we introduce the notation $A(P,F)$ for the affine subspace $\langle P,F\rangle\backslash F$, where $P$ is a point and $F$ and arbitrary subspace. Note that we have $\dim(A(P,F))=\dim(F)+1$. Next, we observe that $(r+1)$-cylinders can be easily constructed starting from a multiset of $q$ points. 

\begin{Construction}
	\label{construction_cylinder}
	Let $r$ and $v'$ be non-negative integers, and consider a $v'$-space $V'$ and a disjoint $r$-space $F$ in $V=\F_q^{v'+r}$. If $\cM' = \{P_1,\dots,P_q\}$ is a multiset of $q$ points in $V'$, then a $(r+1)$-cylinder can be constructed as the multiset $\cM$ consisting of the points of $A(P_i,F)$, $i = 1,\dots,q$.
\end{Construction}

\begin{Proposition}
	The multiset of points of an $(r+1)$-cylinder is $q^r$-divisible.    
\end{Proposition}
\begin{Proof}
	We use the notation of \Cref{def_cylinder} for a given $(r+1)$-cylinder. The statement is trivial for $r=0$ so that we assume $r\ge 1$. Each hyperplane $H$ 
	intersects $F$ either in dimension $r$ or $r-1$. In the first case we have $|(L_i\setminus F) \cap H|\in\left\{0,q^{r}\right\}$. In the second case 
	we have $|(L_i\setminus F) \cap H|=q^{r-1}$ for all $1\le i\le q$. Thus, we have $|\cM\cap H| \equiv 0\pmod{q^{r}}$ for the corresponding 
	multiset of points $\cM$ of the $(r+1)$-cylinder.
\end{Proof}

So, $(r+1)$-cylinders yield $q^r$-divisible multisets of $q^{r+1}$ points and the question arises if there are other isomorphism types. Indeed there are. Any multiset of $q$ (possibly equal) points with multiplicity $q^{r}$ each is $q^{r}$-divisible. For that reason we will consider sets of points instead of multisets in the remaining part. 
We remark that studying multisets of points with restricted point multiplicity might be an interesting problem, but we will not go into this here.  
It will also depend on the dimension whether other isomorphism types exist. Since each set $\cS$ of points in $\aspace$ can be embedded in $\F_q^{v'}$ for 
$v'>v$ we will always assume that $\cS$ is spanning, i.e. spans $\aspace$.

\begin{Observation} 
	In \Cref{construction_cylinder}, $\cM$ is spanning if and only if $\cM'$ is spanning and $\cM$ is a set if and only if $\cM'$ is a set.
\end{Observation}

We call an $(r+1)$-cylinder spanning or projective if the corresponding multiset of points is.

\begin{Lemma} \label{lemma_dimension_cylinder}
	There exists a spanning projective $(r+1)$-cylinder in $\PG(v-1,q)$ if and only if $r+2\le v\le r+q$.
\end{Lemma}
\begin{Proof}
	Due to the above observations it suffices to remark that $2\le \dim(\left\langle \cM'\right\rangle)\le q$ and all dimensions in that range can indeed be attained. 
\end{Proof}

Let $(v,r,q)$ be a triple of non-negative integers where $q$ is a prime power. The question of classification thus boils down to the following.

\begin{Question}\label{conj_generalized_cc}
	If $\cS$ is a $q^r$-divisible spanning set of $q^{r+1}$ points in $\aspace$, does it follow that $\cS$ must be a $(r+1)$-cylinder?
\end{Question}

We say that the generalized cylinder conjecture is true (or false) for a triple $(v,r,q)$ if the answer to this question is affirmative (or not). Notice that the divisibility assumption is as strong as possible and for good reason. It's not hard to construct a $q^r$-divisible set of $q^{n}$ points in $\aspace$, $r+1 < n \leq v-2$, as the disjoint union of $q^{n-r-1}$ different $(r+1)$-cylinders, while not being a cylinder itself.

As we will see in the next section, it makes no difference if we consider point sets in affine or projective geometries. Therefore, the original cylinder conjecture asserts that the answer is affirmative for the triple $(4,1,p)$, $p$ prime. Remark that when $r = 0$, the conjecture is trivially true, so we assume $r \geq 1$ from now on. In \Cref{sec_general_results} we will prove some structural results and generalize ideas from \cite{de2019cylinder}, using a combination of linear programming and geometrical methods. We discuss some cases where we can determine the validity of the generalized cylinder conjecture in \Cref{sec_posneg_results}. Our main results in these two sections can be summarised as follows. In \Cref{lemma_heritable} we have seen that the restriction $\cS \cap H$ of a $q^r$-divisible set of points $\cS$ in $\PG(v-1,q)$ to a hyperplane $H$ is $q^{r-1}$-divisible. We will show that if $|\cS|=q^{r+1}$, then there are many hyperplanes $H$ with $|\cS \cap H|=q^r$, and hence inductive arguments can be applied. Indeed, it will turn out that everything boils down to the case $r = 1$, as \Cref{cor_gcc_equivalence} shows that the generalized cylinder conjecture is true for $(v,r,q)$ if and only if it is true for $(v-r+1,1,q)$. The generalized cylinder conjecture is not always true: if $q$ is a proper prime power, i.e., not a prime, we can construct counterexamples for suitable dimensions based on the existence of subfields, see \Cref{cor_counterexample_nonprime}. For small values of $q$ we find that the generalized cylinder conjecture is true whenever $q\in\{2,3,5\}$, for the triples $(3,1,4)$ and $(4,1,4)$, but not for $(5,1,4)$. The special case $(4,1,q)$ is treated in detail in \Cref{sec_cylinder_conjecture}, where the cases $q=5$ (fixing a flawed proof in \cite{de2019cylinder}) and $q=7$ are fully resolved. Although our numerical data is still rather limited, we state:

\begin{Conjecture}
	The generalized cylinder conjecture is true for $(r+3,r,q)$ and for $(v,r,p)$ if $p$ is a prime.
\end{Conjecture}

Finally, we have to remark that the {\lq\lq}right{\rq\rq} generalization of the cylinder conjecture to non-prime field sizes is a bit unclear in the context of directions determined by a point set, but \Cref{def_cylinder} and \Cref{conj_generalized_cc} are at least reasonable in the context of $q^r$-divisible sets.  

\section{General results and the reduction theorem}\label{sec_general_results}

The linear programming approach is based on the following three linear equations, commonly referred to as the standard equations.

\begin{Lemma}\label{lemma_standard_equations_q}
	
	(See e.g.\ \cite[Lemma 6]{partial_spreads_and_vector_space_partitions})
	
	Let $\cS$ be a set of points in $\aspace$ with $|{\cS}|=n$, and let $a_i$ be the number of hyperplanes in $\aspace$ containing exactly $i$ points of $\cS$ ($0\le i\le n$). Then we have
	\begin{eqnarray}
		\sum_{i=0}^{n}a_i &=& \frac{q^v-1}{q-1},\label{eq_ste1}\\
		\sum_{i=1}^{n}ia_i &=& n\cdot \frac{q^{v-1}-1}{q-1},\label{eq_ste2}\\
		\sum_{i=2}^{n}{i\choose 2} a_i &=& {n\choose 2}\cdot \frac{q^{v-2}-1}{q-1}\label{eq_ste3}.
	\end{eqnarray}
\end{Lemma}

\begin{Proof}
	Double-count incidences of the tuples $(H)$, $(P_1,H)$, and
	$(\{P_1,P_2\},H)$, where $H$ is a hyperplane and $P_1\neq P_2$ are
	points contained in $H$.
\end{Proof}

The set $\{a_i\}_i$ is called the spectrum of $\cS$. The set $\cS$ being spanning is equivalent to $a_n=0$, i.e., no hyperplane contains all points. In that case, the standard equations are equivalent to the first three MacWilliams identities for projective linear codes.

We can adapt \Cref{lemma_standard_equations_q} to our situation of $q^r$-divisible sets of points.

\begin{Lemma}\label{lemma_standard_equations_divisible}
	Let $\cS$ be a $q^r$-divisible spanning set of $q^{r+1}$ points in $\aspace$ and let $a_{iq^r}$ be the number of hyperplanes in $\aspace$
	containing exactly $iq^r$ points of $\cS$ ($0\le i\le q$). Then we
	have
	\begin{eqnarray}
		(q-1)\sum_{i=0}^{q-1}a_{iq^r} &=& q^v-1,\label{eq_sted1}\\
		(q-1)\sum_{i=0}^{q-1}ia_{iq^r} &=& q(q^{v-1}-1),\label{eq_sted2}\\
		(q-1)\sum_{i=0}^{q-1}i(iq^r-1) a_{iq^r} &=& q(q^{r+1}-1)(q^{v-2}-1)\label{eq_sted3}.
	\end{eqnarray}
\end{Lemma}

\begin{Proof}
	We use the equations from \Cref{lemma_standard_equations_q}. Multiplying them by $q-1$, using $n=q^{r+1}$, and taking $q^r$-divisibility into account gives
	\begin{eqnarray}
		(q-1)\sum_{i=0}^{q}a_{iq^r} &=& q^{v}-1,\label{eq_step1}\\
		(q-1)\sum_{i=0}^{q}iq^r a_{iq^r} &=& q^{r+1} \left(q^{v-1}-1\right),\label{eq_step2}\\
		(q-1)\sum_{i=0}^{q}{iq^r\choose 2} a_{iq^r} &=& {q^{r+1}\choose 2} (q^{v-2}-1)\label{eq_step3}.
	\end{eqnarray}
	Finally, dividing \eqref{eq_step2} by $q^r$, \eqref{eq_step3} by $q^r/2$, and recalling $a_{q^{r+1}}=0$ gives the stated result.
\end{Proof}

\begin{Lemma}
	\label{lemma_hyperplane_majority}
	Let $\cS$ be a $q^r$-divisible spanning set of $q^{r+1}$ points in $\aspace$. Then, the number $a_{q^r}$ of 
	hyperplanes with the smallest non-zero number of points is at least $\frac{q^v-1}{q-1}-\left(q^{v-r-1}-q+1\right)$.
\end{Lemma}
\begin{Proof}
	Using the notation from \Cref{lemma_standard_equations_divisible}, Equation~\eqref{eq_sted3} minus $2q^r-1$ times Equation~\eqref{eq_sted2} gives
	$$
	(q-1)\sum_{i=0}^{q-1} q^r i(i-2)a_{iq^r}=-q^r\left((q^v - 1) - (q-1)q^{v-r-1} +(q-1)^2\right).
	$$
	Since $i(i-2)\ge 0$ and $a_{iq^r}\ge 0$ for all $2\le i\le q-1$, we conclude $(q-1)a_{q^r}\ge q^v-1 - (q-1)(q^{v-r-1}-q+1)$. 
\end{Proof}

In other words, almost all hyperplanes contain exactly $q^r$ points. For these hyperplanes we might apply induction and 
assume that they are $r$-cylinders, as we will do later on in the proof of \Cref{thm_increase_r}.

The general idea behind the proof of \Cref{lemma_hyperplane_majority} is the linear programming method based on the standard equations, which is a common technique in finite geometry. To be precise, we maximize or minimize a certain $a_{j}$ under the constraints of \Cref{lemma_standard_equations_divisible}, where we assume that all $a_i\in\mathbb{R}_{\ge 0}$. 
Bounds similar as in \Cref{lemma_hyperplane_majority} for other $a_{iq^r}$ can be obtained easily. 

\begin{Lemma}
	\label{lemma_a0_max}
	Let $\cS$ be a $q^r$-divisible spanning set of $q^{r+1}$ points in $\aspace$. Then, the number $a_0$ of 
	empty hyperplanes is at most $(q^{v-r-1}-q+2)/2$.
\end{Lemma}
\begin{Proof}
	From \Cref{lemma_standard_equations_divisible}, $2q^r$ times Equation~\eqref{eq_sted1} minus $3q^r-1$ times Equation~\eqref{eq_sted2} plus Equation~\eqref{eq_sted3} gives
	$$
	(q-1)\sum_{i=0}^{q-1} q^r(i-1)(i-2)a_{iq^r}=(q-1)(q^{v-1}-q^{r+1}+2q^r).
	$$
	Since $(i-1)(i-2)\ge 0$ and $a_{iq^r}\ge 0$ for all $0\le i\le q-1$, we conclude $2q^ra_0\le (q^{v-1}-q^{r+1}+2q^r)$. 
\end{Proof}

\begin{Lemma}
  \label{lemma_lower_bound_empty_hyperplanes}
  Let $\cS$ be a $q^r$-divisible spanning set of $q^{r+1}$ points in $\aspace$. Then, the number $a_0$ of 
  empty hyperplanes is at least $\frac{q^{v-r-1}-1}{q-1}$.
\end{Lemma}
\begin{Proof}
  Applying \Cref{lemma_standard_equations_divisible}, $q^r(q-1)$ times Equation~\eqref{eq_sted1} minus $q^{r+1}-1$ times 
  Equation~\eqref{eq_sted2} plus Equation~\eqref{eq_sted3} gives
  $$
    (q-1)\sum_{i=0}^{q-1} q^r(i-1)(i-q+1)a_{iq^r}=(q^{v-1}-q^r)(q-1).
  $$
  Since $(i-1)(i-q+1)\le 0$ and $a_{iq^r}\ge 0$ for all $1\le i\le q-1$, we conclude $(q-1) a_0 \ge q^{v-r-1}-1$. 
\end{Proof}

\begin{Corollary}
	\label{cor_affine_space}
	For every $q^r$-divisible spanning set of $q^{r+1}$ points in $\PG(v-1,q)$ there is at least one empty hyperplane.
\end{Corollary}
\begin{Proof}
	Since $\PG(v-1,q)$ contains ${v \brack 1}_q=q^{v-1}+q^{v-2}+\dots+1$ points, we have $v\ge r+2$, so that \Cref{lemma_lower_bound_empty_hyperplanes} 
	gives the stated result.
\end{Proof}

Thus, it makes no difference if we speak about point sets in $\AG(v-1,q)$ or $\aspace$. Another implication of \Cref{cor_affine_space} is that the generalized cylinder conjecture is trivially true when the dimension $v$ is small.

\begin{Proposition}
  \label{prop_dim_at_least_r_plus_3}
  Let $\cS$ be a $q^r$-divisible spanning set of $q^{r+1}$ points in $\aspace$. If $v \leq r+2$, then $v = r+2$ and $\cS\simeq \AG(v-1,q)$.  
\end{Proposition}
\begin{Proof}
  As in the proof of \Cref{cor_affine_space} we conclude $v\ge r+2$, so that $v=r+2$. A single empty hyperplane leaves only $q^{r+1}$ possible points, which all have be to contained in $\cS$.
\end{Proof}

We see that the generalized cylinder conjecture is trivially true for all $(v,r,q)$, where $v\le r+2$. In other words, the classification of $q^r$-divisible spanning sets of $q^{r+1}$ points in 
$\aspace$ is \textit{challenging} for $v\ge r+3$ only. 

With these auxiliary results in hand, we can prove our reduction theorem which essentially states that the validity of the generalized cylinder conjecture depends on the difference $v-r$ and not on the values of $v$ and $r$ itself. As a first step, we show that we can decrease $v$ and $r$ simultaneously and preserve the truthfulness. 

\begin{Proposition}
	\label{prop_cylinder_conjecture_plus_one}
	If the generalized cylinder conjecture is true for $(v+1,r+1,q)$, then it is true for $(v,r,q)$. 
\end{Proposition} 
\begin{Proof}
	If the generalized cylinder conjecture is false for $(v,r,q)$, we can apply the following construction to the corresponding counterexample and obtain a counterexample for $(v+1,r+1,q)$. 
	
	Consider a $q^r$-divisible set $\cS$ of $q^{r+1}$ points in $V:=\aspace$. For a point $P$ outside of the ambient space we 
	consider the new ambient space $V':=\langle V, P\rangle$ and set $\cS'=\left\{A(S,P)\,:\, S\in\cS\right\}$. 
	By construction $\dim(V')=v+1$ and $\cS'$ is a set of $q^{r+2}$ points in $V'$. Now let $H'$ by a hyperplane of $V'$. Either $H'=V$ or $H:=H'\cap V$ 
	is a hyperplane of $V$. In the first case the have $\cS'\cap H'=\cS$, which is of cardinality $q^{r+1}$. In the second case we have 
	$|\cS\cap H|\equiv 0\pmod {q^{r}}$. If $P\le H'$, then $|\cS'\cap H'|=q\cdot |\cS\cap H|\equiv 
	0\pmod {q^{r+1}}$. If $P$ is not contained in $H'$ then each of the $q^{r+1}$ affine lines $A(S,P)$ is met by $H'$ in a 
	single point not equal to $P$, so that $|\cS'\cap H'|=q^{r+1}$. Thus, $\cS'$ is $q^{r+1}$-divisible and one can see that it is not a cylinder, if $\cS$ is not.
\end{Proof}

\begin{Theorem}	\label{thm_increase_r}	
	If the generalized cylinder conjecture is true for $(v,1,q)$, then it is true for all $(v+r-1,r,q)$.
\end{Theorem} 
\begin{Proof}
	Due to \Cref{prop_dim_at_least_r_plus_3} we can assume $v\ge 4$.  We will prove the result by induction on $r$. So, assume that the generalized cylinder conjecture is true for 
	$(v+r-2,r-1,q)$. Let $\cS$ be a spanning $q^r$-divisible  set of $q^{r+1}$ points in $\PG(v'-1,q)$, where $v' = v+r-1$ and $r \geq 2$. Now let $\cF$ be the set of 
	points $F$ such that there exists a point $S\in\cS$ with $A(S,F)\subseteq \cS$. 
		
	We will structure our proof into some intermediate results: \\
	
	\begin{enumerate}
		\item[(1)] For each point $S\in\cS$ there exists an $(r-1)$-space $B$ such that $A(S,B)\subseteq \cS$.
		\item[(2)] $\dim(\langle\cF\rangle)\in\{r-1,r\}$.
		\item[(3)] For each $S_1,S_2\in\cS$ there exists an $(r-1)$-space $B$ such that $A(S_1,B)$ and $A(S_2,B)$ are both contained in $\cS$.
		\item[(4)] Let $S\in\cS$ be a point such that there exist $(r-1)$-spaces $B_1\neq B_2$ with $A(S,B_1),A(S,B_2)\subseteq \cS$. Then, we 
		have $A(S,\langle\cF\rangle)\subseteq\cS$.
		\item[(5)] $\cS$ is an $(r+1)$-cylinder. \\           
	\end{enumerate} 
	For (1) we use \Cref{lemma_hyperplane_majority} to conclude that there are at most $q^{v'-r-1}-q+1$ hyperplanes that do not contain exactly $q^r$ points from 
	$\cS$. Thus, for each point $S\in\cS$ at least one of the $\frac{q^{v'-1}-1}{q-1}$ hyperplanes containing $S$ contains exactly $q^r$ points from $\cS$. 
	If $H$ is such a hyperplane, we can apply induction for $H$ and conclude that $\cS\cap H$ can be partitioned into 
	$\cup_{i=1}^q A\!\left(S_i,B\right)$ for some $(r-1)$-space $B$ and $q$ points $S_i\in\cS$. Of course there exists an index $1\le j\le q$ with $S\in A\!\left(S_j,B\right)$ and hence 
	$A(S,B)=A\!\left(S_j,B\right)$.   
	
	For (2) we note that every point $F\in\cF$ is contained in every empty hyperplane, as $F$ is the only point in $A(S,F)$ not in $\cS$, so that $\langle\cF\rangle$ is contained in the 
	intersection of all empty hyperplanes. Since there are exactly $\frac{q^{v'-\dim(\langle\cF\rangle)}-1}{q-1}$ hyperplanes containing $\langle\cF\rangle$ we conclude 
	from \Cref{lemma_lower_bound_empty_hyperplanes} that $\dim(\langle\cF\rangle)\le r+1$ and moreover if equality holds then every hyperplane containing $\langle\cF\rangle$ is empty. 
	Since $\cS\neq\emptyset$, this is only possible if $\langle\cF\rangle$ is a hyperplane itself, i.e., $v'=r+2$ and $v=3$. Thus, we have $\dim(\langle\cF\rangle)\le r$. 
	If $B$ is a subspace according to (1), then $B\subseteq \cF$ so that $\dim(\langle\cF\rangle)\ge r-1$. 
	
	For (3) we can apply the same idea as in (1). Consider the line $L=\left\langle S_1,S_2\right\rangle$, then $L$ is contained in $\frac{q^{v'-2}-1}{q-1} > q^{v'-r-1}-q+1$ 
	hyperplanes (using $r\ge 2$). So, we can use \Cref{lemma_hyperplane_majority} to conclude the existence of a hyperplane $H$ with $L\le H$ and $|\cS\cap H|=q^r$. Induction on this hyperplane then gives the existence of an $(r-1)$-space $B$ such that $A(S,B)\subseteq \cS\cap H$ for all $S\in\cS\cap H$.      
	
	For (4) we note that $B_1,B_2\subseteq \cF$ implies $\dim(\langle\cF\rangle)\ge r$, so that (2) gives $\dim(\langle\cF\rangle)=r$.      
	
	If $r=2$, then take any point $F$ on the line $\langle\cF\rangle$ and a point $S \in \cS$. We will directly prove that $A(S,F) \subseteq \cS$. We know that there are at 
	most $q^{v'-3}-q+1$ hyperplanes not intersecting $\cS$ in $q^2$ points by \Cref{lemma_hyperplane_majority}, so out of the $q^{v'-3}$ hyperplanes through $\langle S,F \rangle$ 
	intersecting $\langle\cF\rangle$ in only $F$, there must be at least $q-1$ hyperplanes containing exactly $q^2$ points. So take one such hyperplane, apply induction and find that 
	$A(S,F)\subseteq \cS$. It follows that $A(S,\langle\cF\rangle) \subseteq \cS$ for all $S \in \cS$, i.e., (4) is valid for $r=2$.
	
	Now assume $r\ge 3$. We can find a point $S \in \cS$ and distinct $(r-1)$-spaces $B_1,B_2$ such that $A(S,B_1)$ and $A(S,B_2)$ both are contained in $\cS$. Now, take a point 
	$F \in \langle\cF\rangle \setminus \{B_1,B_2\}$, and consider any $3$-space $\pi$ through $\langle S,F \rangle$ not intersecting $B_1 \cap B_2$. Then this $3$-space $\pi$ 
	intersects $B_1$ and $B_2$ each in a point, say $F_1$ and $F_2$. There are $\frac{q^{v'-3}-1}{q-1}$ hyperplanes through this $3$-space, which is more than $q^{v'-r-1}-q+1$ 
	if $r \ge 3$. We again conclude by induction that there must be a $(r-1)$-space $B$ such that $A(S,B) \subseteq \cS$. As $F_1$ and $F_2$ must be contained in $B$, we conclude 
	that $F$ must also be and hence $A(S,F) \subseteq \cS$. It again follows that $A(S,\langle\cF\rangle) \subseteq \cS$ for all $S \in \cS$ and so (4) is valid for $r \geq 3$ too.
	
	For the final step (5) we can assume that there exists a point $S\in \cS$ such that there is a unique $(r-1)$-space $B$ satisfying $A(S,B)\subseteq \cS$, as  otherwise by (4) we can already conclude that $\cS$ is an $(r+1)$-cylinder. Under this assumption, it follows from (3) that $A(S',B)\subseteq\cS$ for all $S'\in\cS$, so that modulo $B$ we obtain a set $\cS'$ of $q^2$ points that is $q$-divisible and spans a space of dimension $v'-r+1=v$. For this set $\cS'$ we can apply the generalized cylinder conjecture for $(v,1,q)$ and conclude $\cS'= \cup_{i=1}^q A(S'_i,B')$ for some points $S'_i$ and $B'$. By construction, we then have $\cS=\cup_{i=1}^q A(S'_i,\langle B,B'\rangle)$ which shows that $\cS$ is an $(r+1)$-cylinder.
\end{Proof} 

Combining the previous results, we find the promised reduction theorem.

\begin{Corollary}
	\label{cor_gcc_equivalence}
	The generalized cylinder conjecture is true for $(v,r,q)$ if and only if it is true for $(v-r+1,1,q)$.
\end{Corollary}

Finally, we can transfer many of the insights of \cite{de2019cylinder} to our more general situation.

\begin{Lemma}
  \label{lemma_1}
  ({\em Cf.} \cite[Lemma 1]{de2019cylinder})\\
  Let $\cS$ be a $q^r$-divisible spanning set of $q^{r+1}$ points in $\aspace$. Let $K$ be a subspace of codimension $2$ in $\aspace$. Assume that $|\cS\cap K|=kq^{r-1}$ for some 
  integer $0<k<q$. Then, any hyperplane containing $K$ contains at most $kq^r$ points from $\cS$.
\end{Lemma}
\begin{Proof}
  Since every hyperplane containing $K$ contains at least $kq^{r-1}$ points, it should contain at least $q^r$ points. Therefore, counting the number of points on 
  hyperplanes containing $K$, we find at least $(q+1)\left(q^r-kq^{r-1}\right)+kq^{r-1}=q^{r+1}-(k-1)q^r$ points. Hence, there are $(k-1)q^r$ points left, which implies that 
  a single hyperplane contains at most $kq^r$ points.
\end{Proof}

\begin{Corollary}
  \label{corollary_1}
  ({\em Cf.} \cite[Corollary 1]{de2019cylinder})\\
  Let $\cS$ be a $q^r$-divisible spanning set of $q^{r+1}$ points in $\aspace$. Suppose the hyperplane $H$ contains $kq^r$ points of $\cS$, where $0<k<q$, then 
  every hyperplane $K$ of $H$, i.e., $K\le H$ is a subspace of codimension $2$ in $\aspace$, contains either $0$ or $lq^{r-1}$ points of $\cS$ for some $k \leq l \leq q$.
\end{Corollary}
\begin{Proof}
  By \Cref{lemma_heritable} we know that $|\cS\cap K|\equiv 0\pmod{q^{r-1}}$ holds . Suppose that $K$ contains $lq^{r-1}$ points 
  for some integer $0<l<k$. Then, by \Cref{lemma_1}, every hyperplane $H$ in $\aspace$ containing $K$ contains at most $lq^r<kq^r$ points, which 
  is a contradiction.  
\end{Proof}

For the next proof, denote by $[x]_q := (q^x-1)/(q-1)$ the number of hyperplanes through a codimension $x$-space.

\begin{Theorem}{({\em Cf}.~\cite[Theorem 2]{de2019cylinder})}\\
	\label{thm_cy_full_space}
	Let $\cS$ be a $q^r$-divisible spanning set of $q^{r+1}$ points in $\aspace$. If $\cS$ contains a full affine $(r+1)$-space then $\cS$ is a $(r+1)$-cylinder.
\end{Theorem}
\begin{Proof}
	Due to \Cref{prop_dim_at_least_r_plus_3} we can assume $v \geq r+3$. Denote by $A$ the affine $(r+1)$-space contained in $\cS$ and let $A_\infty$ be `its part at infinity', i.e. the unique $r$-space so that $A \cup A_\infty$ is an $(r+1)$-space of $\aspace$. Take a point $P \in \cS \setminus A$ and let $H$ be a hyperplane containing $A$ and $P$. By $q^r$-divisibility, we can assume that $|\cS \cap H|=kq^r$ with $1<k<q$. 
	
	Fix a point $Q \in A$ and consider the set $\cK$ of $(v-2)$-spaces through $Q$, not containing $A_\infty$. Observe that $|\cK| = [v-2]_q - [v-r-2]_q$ and that every point in $(\cS \cap H) \setminus A$ is contained in $a := [v-3]_q - [v-r-3]_q$ of these spaces. Now double counting pairs $\{(K,R) \,\, | \,\, K \in \cK, R \in K \cap \cS \setminus A\}$ we find
	\[(kq^r-q^r)a = \sum_{K \in \cK}|(K \cap \cS \setminus A)| \geq |\cK|(kq^{r-1}-q^{r-1}),\]
	by \Cref{corollary_1}. As $|\cK| = aq$, we see that in fact we have equality so that $|K \cap \cS| = kq^{r-1}$ for all $K \in \cK$. Moreover, the point $Q \in A$ was arbitrary so we conclude that every $(v-2)$-space not through $A_\infty$ contains $kq^{r-1}$ points.
	
	Now retaking the point $P \in \cS \cap H \setminus A$, we can consider the set $\cK'$ of $(v-2)$-spaces through $P$, not containing $A_\infty$. Again, this set has size $|\cK'| = [v-2]_q - [v-r-2]_q$. Denote by $B$ the $(r+1)$-space $\langle A_\infty,P \rangle$ and $x = |\cS \cap B|$. If we can show that $x = q^r$, we are done. We can now double count the pairs $\{(K,R) \,\, | \,\, K \in \cK', R \in K \cap \cS \setminus \{P\}\}$. If $R \in B$ then the number of $(v-2)$-spaces through $P$ and $R$, not containing $A_\infty$ is $b = [v-3]_q-[v-r-2]_q$. For a point not in $B$, the correct count for the $(v-2)$-spaces is still $a$ and we find the following equality: 
	
	\[(kq^r-x)a + (x-1)b = \sum_{K \in \cK'}|(K \cap \cS \setminus B)| = |\cK'|(kq^{r-1}-1).\]
	
	After simplification we indeed find $x = q^r$, concluding the proof.
	
%	\[(kq^r-x)b + (x-1)c = \sum_{K \in \cK'}|(K \cap \cS \setminus B)| = bq(kq^{r-1}-1),\]
%	
%	\[(b-c)x + c = bq\]
%	
%	as $b-c = q^{v-r-3}$ and $bq - c = q^{v-3}$ we find $x = q^r$.
\end{Proof}

\section{Positive and regative results}\label{sec_posneg_results}

In this section we will discuss some triples $(v,r,q)$ for which we can determine the validity of the generalized cylinder conjecture. Before that, we introduce some notation to make the proofs less cumbersome. 

If $K$ is a $k$-space in $\aspace$, then $K(\cS) := |K \cap \cS|$ denotes the number of points in the intersection, which we will also refer to as the multiplicty of $K$. If $K$ is a point, then we will say that $K$ is a $0$- or $1$-point, whenever $K(\cS) = 0$ or $1$ respectively. Similarly, if $K$ is a line or a plane, we will say that $K$ is an $m$-line or $m$-plane if $K(\cS) = m$. 

First of, we will show that the generalized cylinder conjecture is true whenever $q\in\{2,3\}$.

\begin{Lemma}
	\label{lemma_gcc_two_weights}
	%%For some non-negative integer $r$ 
	Let $\cS$ be a spanning set of $q^{r+1}$ points in $\PG(v-1,q)$. If every hyperplane of $\PG(v-1,q)$ contains 
	either $q^r$ or no point from $\cS$, then $v=r+2$ and $\cS\simeq \AG(v-1,q)$.  
\end{Lemma}
\begin{Proof}
	Since only $a_0$ and $a_{q^r}$ are non-zero in \Cref{lemma_standard_equations_divisible}, we immediately find $v=r+2$, so that we can apply \Cref{prop_dim_at_least_r_plus_3}.
\end{Proof}

\begin{Corollary}
	\label{cor_gcc_q_2}
	The generalized cylinder conjecture is true for all triples $(v,r,2)$.
\end{Corollary}	

Next up, we will show the generalized cylinder conjecture for $(v,1,3)$. By \Cref{cor_gcc_equivalence}, this suffices to conclude that the generalized cylinder conjecture is true for all triples $(v,r,3)$. 

\begin{Proposition}
	\label{prop_gcc_q_3}
	Let $\cS$ be a $3$-divisible spanning set of $9$ points in $\PG(v-1,3)$. Then $\cS$ is a $2$-cylinder.
\end{Proposition}
\begin{Proof}
	Assume to the contrary, that $\cS$ is not a $2$-cylinder, so that \Cref{prop_dim_at_least_r_plus_3} implies $v\ge 4$. Since the maximum multiplicity 
	of a hyperplane is $6$, each subspace of multiplicity $6$ is a hyperplane. Assume that $K$ is a subspace of multiplicity $4$, so that $\dim(K)\le v-2$. 
	We denote the codimension $v-\dim(K)$ of $K$ by $x$, and hence $x\ge 2$. Since there are $[x]_3$ hyperplanes through $K$, every $1$-point outside of $K$ 
	is contained in $[x-1]_3$ hyperplanes, and every hyperplane through $K$ has multiplicity $6$, we have $2[x]_3=5[x-1]_3$, so that $3^{x-1}=-3$, which is impossible. 
	Thus, no subspace can have a multiplicity of exactly $4$. 
	
	By \Cref{lemma_gcc_two_weights} we can assume the existence of 
	%From the standard equations we compute $a_0=\left(3^{v-2}-1\right)/2$, 
	%$a_3=\left(7\cdot 3^{v-2}+3\right)/2$, and $a_6=\left(3^{v-2}-3\right)/2$, so that $a_6\ge 3$.
	a hyperplane $H$ with $\cS(H)=6$. Since $\cS$ is spanning we have $\dim(H)\le 6$. If $\dim(H)=6$, then we can assume w.l.o.g.\ that the $1$-points in $H$ are given by $\left\langle e_1\right\rangle, \dots, \left\langle e_6\right\rangle$, 
	where the $e_i$ denote the standard unit vectors. With this, the subspace $\left\langle e_1,\dots,e_4\right\rangle$ would have multiplicity $4$, which is a contradiction. Now 
	assume $\dim(H)=5$ and that the $1$-points in $H$ are given by $\left\langle e_1\right\rangle, \dots, \left\langle e_5\right\rangle$ and a sixth point $P$. Consider the 
	subspace $\left\langle e_1,\dots,e_4\right\rangle$. Since it does not contain $e_5$ and there is no subspace of multiplicity four, it has to contain $P$, which means that the fifth coordinate of the vectors in $P$ are zero. We can repeat this argument for the subspace $\left\langle \left\{e_1,\dots,e_5\right\}\backslash e_i\right\rangle$, $i = 1,\dots,4$ and conclude that the $i$-th coordinate of the vectors in $P$ are zero for all $1\le i\le 5$, which is a contradiction. Thus, the remaining possibilities are 
	$\dim(H)\in\{3,4\}$, i.e., $v\in\{4,5\}$. From \Cref{thm_cy_full_space} we conclude that the maximum line multiplicity is at most $2$. So, if 
	$\dim(H)=3$, then $\cS \cap H$ would be a set of $6$ points in $\PG(2,3)$ with line multiplicity at most two, which does not exist. Therefore, we have $\dim(H)=4$ and $v=5$. Using \Cref{corollary_1} we 
	conclude that every plane $\pi$ in $H$ has a multiplicity in $\{0,2,3,5\}$. By the same reasoning as before, we cannot have $5$ points in $\PG(2,3)$ with line multiplicity at most two and hence $\cS(\pi) \neq 5$. For the 
	spectrum $\left(a_i'\right)$ of $\cS \cap H$ the standard equations yield the unique solution $a_0' = 8$, $a_2' = 18$, $a_3' = 14$. Now consider the subspaces spanned by 
	one of the ${6\choose 3}=20$ triples of $1$-points in $H$. As no line contains $3$ points and every plane contains at most $3$ points, each triple of points spans a distinct plane, implying $a_3'\ge 20$, which is a contradiction.           
\end{Proof}

\begin{Corollary}
	\label{cor_gcc_q_3}
	The generalized cylinder conjecture is true for all triples $(v,r,3)$.
\end{Corollary}

An interesting implication is that the dimension $v$ of every $3^r$-divisible projective $\left[3^{r+1},v\right]_3$-code is at most $r+3$ for every positive integer $r$. 
We remark that Ward's upper bound on the dimension of divisible codes \cite[Theorem 6]{ward2001divisible} gives $v \leq q(r+1)$ for the triple $(v,r,q)$, and is hence not strong enough to give this result. Using the software package \texttt{LinCode} \cite{kurz2019lincode} we have computationally checked that there is no $3$-divisible $[9,\ge 5]_3$-code, 
no $9$-divisible $[27,\ge 6]_3$-code, and no $27$-divisible $[81,\ge 7]_3$-code. This means it might not be necessary to assume that $\cS$ is a set and not a multiset to obtain the stated upper bound for the dimension. We remark that the truth of the cylinder conjecture for $(4,1,2)$ and $(4,1,3)$ was also proven in \cite{de2019cylinder}. 

In principle it is possible to enumerate all projective $q^r$-divisible $\left[q^{r+1},v\right]_q$ codes and to check whether the corresponding point sets are $(r+1)$-cylinders 
for given finite parameters. However, given the currently available software for the exhaustive enumeration of linear codes, this approach is limited to rather small parameters. 
Nevertheless we report our corresponding findings here. The last step -- checking whether all resulting point sets are $(r+1)$-cylinders -- can be replaced by a 
counting argument. The numbers of projective linear codes over $\mathbb{F}_5$ of effective lengths $n=5$ ordered by their dimension $k$ are given by $2^1 3^4 4^3 5^1$, as 
can be easily enumerated using the software package \texttt{LinCode} \cite{kurz2019lincode} -- even a classification by hand is possible. So, \Cref{construction_cylinder} 
yields  $3^1 4^4 5^3 6^1$ $5$-divisible projective linear codes over $\mathbb{F}_5$ of effective lengths $n=25$, again ordered by their dimension $k$. Using \texttt{LinCode} 
we verified that there are no further $5$-divisible projective linear codes over $\mathbb{F}_5$ of effective length $n=25$. Thus, we have computationally proven  
that the generalized cylinder conjecture is true for $(v,1,5)$, where the dimension $v$ is arbitrary. This covers the special case $(4,1,5)$ that we treat in the subsequent section.  
From \Cref{cor_gcc_equivalence} we conclude:
\begin{Corollary}
	\label{cor_gcc_q_5}
	The generalized cylinder conjecture is true for all triples $(v,r,5)$.
\end{Corollary}

For $q\in\{2,3,4\}$ we can perform the same computation. The cases $q=2,3$ verify our theoretical findings for $(v,1,q)$. The number of projective linear codes 
over $\F_4$ of effective lengths $n=4$ ordered by their dimension 
$k$ are given by $2^1 3^2 4^1$. The number of $4$-divisible projective linear codes over $\mathbb{F}_4$ of effective lengths $n=16$ ordered by their dimension $k$ are 
given by $3^1 4^2 5^2$. In other words, the generalized cylinder conjecture is true for $(3,1,4)$ and $(4,1,4)$ but not for $(5,1,4)$. For $v\ge 6$ there do not exist projective 
$4$-divisible $[16,v]_4$-codes so that the generalized cylinder conjecture is trivially true for $(v,1,4)$ whenever $v\ge 6$. Therefore, \Cref{cor_gcc_equivalence} gives:
\begin{Corollary}
	\label{cor_gcc_q_4}
	The generalized cylinder conjecture is true for $(v,r,4)$ if and only if $v\neq r+4$.
\end{Corollary}

Of the two linear codes in dimension $5$, the one that does not correspond to a $2$-cylinder has a generator matrix given by
$$
\begin{pmatrix}
	0110101101110000\\
	1101100011101000\\
	1100011111000100\\
	1111111000000010\\
	0111010110100001
\end{pmatrix}.
$$
The code has weight enumerator $W(z)=1z^{0}+90z^{8}+840z^{12}+93z^{16}$ and an automorphism group of order $1935360$. Considered over $\F_2$ the stated generator matrix 
gives a linear $\F_2$ code with weight enumerator $W(z)=1z^{0}+30z^{8}+1z^{16}$, which means that the code is an affine $5$-space.

Computationally we also verified that the generalized cylinder conjecture is true for $(5,2,4)$, which also follows from \Cref{thm_increase_r}. Due to the counter example for $(5,1,4)$ there 
is also a counter example for $(6,2,4)$, see \Cref{prop_cylinder_conjecture_plus_one} . We have computationally checked that this counter example is unique. 

In order to generalize the above counter example to the generalized cylinder conjecture for $(5,1,4)$ we remark that for each integer $h\ge 2$ the field $\F_q$ is a subfield of 
$\F_{q^h}$, so that $\F_{q^h}^v\cong \F_q^{vh}$. Using this isomorphism we can we can embed every multiset of points $\cM'$ in $\F_q^v$ as a multiset of points $\cM$ 
(of the same cardinality) in $\F_{q^h}^v$. Moreover, every $k$-space in $\F_{q^h}^v$ corresponds to a $kh$-space in $\F_q^{vh}$. 

\begin{Lemma}
	\label{lemma_gcc_counter_example}
	Let $\cS'$ be a spanning projective $2h$-cylinder in $\F_q^v$, where $h\ge 2$. Then the corresponding embedding $\cS$ in $\F_{q^{h}}^v\cong \F_q^{vh}$ 
	is a spanning projective $q^h$-divisible set of $\left(q^h\right)^2$ points in $\PG\!\left(v-1,q^h\right)$ that is not a $2$-cylinder  
\end{Lemma}
\begin{Proof}
	First we observe $|\cS|=|\cS'|=q^{2h}=\left(q^h\right)^2$. Since $\cS'$ is spanning and projective, the same applies to $\cS$. 
	An arbitrary hyperplane $H$ in $\F_{q^h}^v$ has dimension $(v-1)h$ over $\F_q$. Since $\F_q^v$ has dimension 
	$v$ over $\F_q$, there exists a subspace $K$ in $\F_q^v$ of dimension at least $v-h$ such that $\left|\cS\cap H\right|=\left|\cS'\cap K\right|$. Note that 
	$\cS'$ is $q^{2h-1}$-divisible, so that $\left|\cS'\cap K\right|\equiv 0\pmod {q^h}$ due to \Cref{lemma_heritable}. 
	and we conclude that $\cS$ is $q^h$-divisible.  
	
	Now assume that $\cS$ is a $2$-cylinder and let $L$ be one of the $q^h$-lines. Consider two $1$-points $P_1,P_2$ on $L$ and denote by $P_1',P_2'$ be the corresponding points in $\cS'$. The line $L'=\left\langle P_1',P_2'\right\rangle$ has multiplicity at most $q$ in $\cS'$, so that $L$ has a multiplicity of at most $q$ in $\cS$, which is a contradiction due to $h\ge 2$.   
\end{Proof}

\begin{Corollary}\label{cor_counterexample_nonprime}
	 For any integer $h\ge 2$, the generalized cylinder conjecture is false for $\left(v,r,q^h\right)$ whenever $2h+r\le v\le 2h+r+q-2$.
\end{Corollary}
\begin{Proof}
	From \Cref{lemma_gcc_counter_example} and \Cref{lemma_dimension_cylinder} we conclude that the generalized cylinder conjecture is wrong for 
	$\left(v,1,q^h\right)$, where $2h+1\le v\le 2h-1+q$, so that \Cref{cor_gcc_equivalence} gives the general statement. 
\end{Proof}

If $v=2h+r$, then the point set $\cS'$ in \Cref{lemma_gcc_counter_example} is an affine geometry, so that $\cS$ is an affine subgeometry. For the 
special case $h=2$ one also speaks of a Baer (sub-)geometry. In general, our construction is an instance of the technique of the so-called 
\emph{field reduction}, which yields a lot of non-trivial constructions and characterizations of geometric and algebraic structures, see 
e.g.~\cite{lavrauw2015fieldreduction}. Of course one might conjecture that every $q^r$-divisible set of $q^{r+1}$ points in $\PG(v-1,q)$ is either an $(r+1)$-cylinder or arises from a cylinder over a subfield.  

\section{The generalized cylinder conjecture for $(4,1,q)$}
\label{sec_cylinder_conjecture}

In this section, we will focus on the case $v = r+3$, and by \Cref{cor_gcc_equivalence} we can restrict ourselves to the triple $(4,1,q)$. We will gather some more information on a possible counterexample, which leads us to be able to prove the generalized cylinder conjecture for $(4,1,5)$ and $(4,1,7)$. A proof of the former was also claimed in \cite{de2019cylinder}, but the proof contains an error which we correct here.

For the results in this section, we will often make the following assumption:

\[(\star) \hspace{1cm} \cS \text{ is a } q\text{-divisible spanning set of } q^2 \text{ points in } \PG(3,q) \text{ which is not a 2-cylinder.}\]

\vspace{1em}

Since $r$ is fixed in this section, we will also refer to a 2-cylinder as a cylinder. Moreover, we can consider $\cS$ as a set of points in $\AG(3,q)$ by \Cref{cor_affine_space}. Lastly, by \Cref{thm_cy_full_space} we can assume that any plane intersects $\cS$ in at most $q^2-q$ points. Our general strategy is to obtain some structural results on $\cS$ and find a contradiction for small $q$. We will heavily rely on the standard equations for planes and points in $\cS$ as stated in \Cref{lemma_standard_equations_divisible}, but also similar equations for lines in a $kq$-plane $H$ and points in $\cS \cap H$, obtained by the same double counting method. The latter gives us information on the number of $i$-lines for each $i$, which we will also refer to as the spectrum.

We start off by investigating the multiplicity of a line with respect to $\cS$. Summarizing the conclusions of \Cref{cor_affine_space}, \Cref{corollary_1} and \Cref{thm_cy_full_space}, we can state the following lemma.

\begin{Lemma}
	\label{lemma_cy_line_multiplicities_in_plane}
	Under $(\star)$, the line multiplicities in a $kq$-plane $H$ are contained in $\{0,k,k+1,\dots,q-1\}$. 
\end{Lemma}

The restriction on the possible line multiplicities is quite severe. Indeed, we can investigate the existence of a set $\cK$ of $kq$ points in $\PG(2,q)$ admitting these line multiplicities, independently of the (generalized) cylinder conjecture. It turns out that when $k$ is large, such sets cannot exist and hence $\cS$ cannot have large intersections with planes. This idea is illustrated in the next few results.

\begin{Lemma}
	\label{lemma_cy_no_qm1_plane}
	For a set of $q(q-1)$ points in $\PG(2,q)$, there exists a line with a multiplicity not in $\{0,q-1\}$.
\end{Lemma}
\begin{Proof}
	Otherwise the first two standard equations would give $a_0=1$ and $a_{q-1}=q^2+q$, so that the third one yields the contradiction
	$$
	{{q-1}\choose 2}a_{q-1}=\frac{q(q-1)(q^2-q-2)}{2}<\frac{q(q-1)(q^2-q-1)}{2}={{q(q-1)}\choose 2}.   
	$$
\end{Proof}

\begin{Corollary}
	\label{cor_cy_no_qm1_plane}
	Under $(\star)$, there can be no $q(q-1)$-plane.
\end{Corollary} 
\begin{Proof}
	If $H$ is a $q(q-1)$-plane, then \Cref{lemma_cy_line_multiplicities_in_plane} and \Cref{lemma_cy_no_qm1_plane} yield a 
	contradiction.   
\end{Proof}

We can use this result to find an alternative proof of the cylinder conjecture when $q \in \{2,3\}$.

\begin{Corollary}
	\label{cor_cy_conj_true_2_3}
	The generalized cylinder conjecture is true for the triples $(4,1,2)$ and $(4,1,3)$.
\end{Corollary}
\begin{Proof}
	Assume that $\cS$ is not a cylinder, so that \Cref{cor_cy_no_qm1_plane} implies $a_{q(q-1)}=0$. For $q=2$ this means that all points 
	of $\cS$ are contained in $0$-planes, which is absurd. For $q=3$ solving the first two standard equations give $a_0=1$ and $a_3=39$. The third 
	implies the contradiction
	$$
	{3\choose 2}a_3=117 <144={9\choose 2}\cdot \left(3+1\right).   
	$$   
\end{Proof}

\begin{Remark}
	The argumentation of the proof of \Cref{cor_cy_conj_true_2_3}, i.e., the first three standard equations combined with \Cref{cor_cy_no_qm1_plane}, would 
	give the unique spectrum $a_0=7$, $a_1=72$, and $a_2=6$ for $q=4$.
\end{Remark}

\begin{Lemma}
	\label{lemma_cy_no_qm2_plane}
	For a set of $q(q-2)$ points in $\PG(2,q)$, $q \geq 4$, there exists a line with a multiplicity not in $\{0,q-2,q-1\}$.
\end{Lemma}
\begin{Proof}
	Otherwise the first two standard equations would give $a_{q-1}=a_0(q-2)-(q-2)$ and $a_{q-2}=q^2+2q-1-a_0(q-1)$,
	so that the third yields
	$$
	0=a_{q-2}{{q-2}\choose 2}+a_{q-1}{{q-1}\choose 2}-{{q(q-2)}\choose 2}=\frac{(q-2)(a_0(q-1)-3q+1)}{2},
	$$ 
	which implies
	$$
	a_0=\frac{3q-1}{q-1}=3+\frac{2}{q-1}\notin\N_0.
	$$
\end{Proof}

\begin{Remark}
	In the proof of \Cref{lemma_cy_no_qm2_plane} the third standard equation is needed since e.g.\ $\left(a_0,a_{q-2},a_{q-1}\right)=(2,17,2)$ satisfies the 
	first two standard equations for $q=4$. 
\end{Remark}

\begin{Corollary}
	\label{cor_cy_no_qm2_plane}
	Under $(\star)$, there can be no $q(q-2)$-plane.
\end{Corollary} 
\begin{Proof}
	Due to \Cref{cor_cy_conj_true_2_3} we can assume $q\ge 4$. If $H$ is a $(q-2)q$-plane, then \Cref{lemma_cy_line_multiplicities_in_plane} and \Cref{lemma_cy_no_qm2_plane} yield a contradiction.   
\end{Proof}

\begin{Corollary}
	\label{cor_cy_conj_true_4}
	The generalized cylinder conjecture is true for the triple $(4,1,4)$.
\end{Corollary}
\begin{Proof}
	Assume that $\cS$ is not a cylinder, so that \Cref{cor_cy_no_qm1_plane} and \Cref{cor_cy_no_qm2_plane} imply 
	$a_{q(q-1)}=0$ and $a_{q(q-2)}=0$. For $q=4$ solving the first two standard equations give $a_0=1$ and $a_4=84$. The third 
	implies the contradiction
	$$
	{4\choose 2}a_4=504 <600={16\choose 2}\cdot \left(4+1\right).   
	$$   
\end{Proof}

We might continue in the vein of \Cref{lemma_cy_no_qm2_plane} and consider sets of $q(q-3)$ points in $\PG(2,q)$ whose line multiplicities are contained in $\{0,q-3,q-2,q-1\}$. 
Due to to \Cref{cor_cy_conj_true_2_3} and \Cref{cor_cy_conj_true_4} we are only interested in the cases where $q\ge 5$. It turns out that he unique possibility is given by 
$q=5$ and a spectrum given by $a_0=6$, $a_2=15$, $a_3=10$, and $a_4=0$. So, similar to \Cref{cor_cy_no_qm2_plane},  
we can conclude that a $q(q-3)$-plane does not exist unless $q=5$. This result also follows in a much shorter way from \Cref{lemma_cy_gamma_2_at_most_q_minus_2}, but we want to remark that the latter is not necessary to show that the generalized cylinder conjecture is true for the triple $(4,1,5)$.

%Solve standard equations in function of $a_0$ and find that $a_0 \geq 4 + 3/(q-1)$ and $a_0 \leq 4 + 6/(q-2)$. So if $q \geq 9$ this is a contradiction, otherwise the only possibilities are $q = 5$ and spectra (5,21,2,3) or (6,15,10,0), or $q = 7$ with spectrum (5,38,12,2), or $q = 8$ with spectrum (5,48,20,0). Now consider a $q-2$ line and suppose $a_0 = 5$, then there is a $0$-point with at most one $0$-line through it. Counting the points on lines through this $0$-point leads to $q(q-3)+1$ points. We conclude that $a_0 = 5$ is not possible.

\begin{Proposition}
	\label{prop_cy_conj_true_5}
	The generalized cylinder conjecture is true for the triple $(4,1,5)$. 
\end{Proposition}
\begin{Proof}
	By \Cref{cor_cy_no_qm1_plane} and \Cref{cor_cy_no_qm2_plane} we know that 
	$a_{q(q-1)}=0$ and $a_{q(q-2)}=0$. With this, the standard equations for $\cS$ yield the unique spectrum $a_0=11$, $a_5=135$, 
	and $a_{10}=10$. Now let $H$ be a $10$-plane and $b_i$ the number of lines of multiplicity $i$, where $b_i=0$ for $i\notin\{0,2,3,4\}$. From the standard 
	equations we conclude $b_{2} = 51-6b_0$, $b_{3} = -38+8b_0$, and $b_{4} = 18 -3b_0$, so that $b_4\ge 0$ implies $b_0\le 6$ and 
	$b_3\ge 0$, $b_0\in\N_0$ imply $b_0\ge \left\lceil \frac{19}{4}\right\rceil =5$. Assume $b_0=5$ and note that two $4$-lines 
	cannot share a common $1$-point $P$ since otherwise counting points on the lines $L_1,\dots,L_{6}$ through $P$ would yield the contradiction
	$$
	10=\cS(H)=\sum_{i=1}^{6} \cS(L_i)-5\cdot\cS(P)\ge 2\cdot 4+4\cdot 2-5=11
	$$
	using $\cS(L_i)\ge 2$ for all $1\le i\le 6$. Thus, the $4$-lines are pairwise disjoint, so that 
	$ 10=\cS(H) \ge  a_{4}\cdot 4= 12$, which is a contradiction. Hence, we have $b_0=6$, $b_2=15$, $b_3=10$, and $b_i=0$ otherwise.
	
	For a line $L$ of multiplicity $3$ in $H$ denote the other $5$ planes by 
	$H_1,\dots,H_5$. Since $\cS(H_i)\in\{5,10\}$ for all $1\le i\le 5$ and $|\cS|=25$, there exists an index $1\le i\le 5$ with $\cS(H_i)=10$. 
	Due to $1+b_3\cdot 1=11$, there are at least eleven $10$-planes, which contradicts $a_{10}=10$. 
\end{Proof}

\begin{Remark}
	\label{remark_10_0_2_3_arc_pg_2_5}
	There exists a unique projective $[10, 3, \{6, 7, 8, 10\}]_5$ code $\cC$ with generator matrix 
	$$
	\begin{pmatrix}
	1&1&1&1&1&1&0&1&0&0\\
	4&4&3&2&1&0&1&0&1&0\\
	3&4&4&2&0&1&4&0&0&1
	\end{pmatrix}.
	$$    
	This code has an automorphism group of order $480$ and weight enumerator $W_{\cC}(x)= 1x^0 +40x^7 +60x^8 +24x^{10}$. A geometrical interpretation for this set of points was also given in \cite{vandevoorde2011}.
	We remark that the existence of the above code was excluded in the proof of \cite[Theorem 4]{de2019cylinder} and so the proof is flawed. More precisely, in the last sentence of the argument showing the existence of a $4$-line in a $10$-plane where two further $0$-lines are constructed, 
	it can happen that they (partially) coincide with the four $0$-lines found before.
\end{Remark}

\begin{Lemma}
	\label{lemma_cy_gamma_2_at_most_q_minus_2}
	Under $(\star)$, the maximum line multiplicity with respect to $\cS$ is $q-2$.
\end{Lemma}
\begin{Proof}
	Due to \Cref{cor_cy_conj_true_2_3} and \Cref{cor_cy_conj_true_4} we can assume $q\ge 5$. 
	By \Cref{lemma_cy_line_multiplicities_in_plane} we already know that the maximum line multiplicity is at most $q-1$.  
	So assume that $L$ is a line of multiplicity $q-1$ and denote by $Q_1,Q_2$ the two $0$-points on $L$. Let $H$ be an arbitrary 
	hyperplane containing $L$, where $\cS(H)=kq$ for an integer $1\le k<q$. 
	
	For each $1$-point $P$ on $L$ let $L_1,\dots,L_q$ denote the $q$ lines trough $P$ in $H$ that are not equal to $L$. Note that 
	\Cref{lemma_cy_line_multiplicities_in_plane} implies $\cS(L_i)\ge k$ for all $1\le i\le q$. From
	$$
	kq=\cS(H)=\cS(L)+\sum_{i=1}^{k} \cS(L_i)-q\cdot\cS(P)\ge q-1+qk-q=kq-1
	$$
	we conclude that $q-1$ of the $L_i$, where $1\le i\le q$, have multiplicity $k$ and one has multiplicity $k+1$.   
	
	Now consider a $0$-point $R$ in $H$ not on $L$ and let $L_i'$ be the lines through $R$ and $Q_i$ in $H$, where $1\le i\le 2$. By $L_3',\dots,L_{q+1}'$ we denote the 
	remaining $q-1$ lines through $R$ in $H$. Note that the $L_i'$ meet the line $L$ in a $1$-points, so that $\cS(L_i')\in\{k,k+1\}$ for 
	$3\le i\le q+1$. From
	$$
	kq=\cS(H)=\sum_{i=1}^{q+1} \cS(L_i') -q\cdot\cS(R) \ge \cS(L_1')+\cS(L_2')+(q-1)k  
	$$ 
	we conclude $\cS(L_1')+\cS(L_2')\le k$, so that $\cS(L_1'),\cS(L_2')\in\{0,k\}$ due to \Cref{lemma_cy_line_multiplicities_in_plane}. Thus, 
	for $1\le i\le 2$ the $q+1$ lines through $Q_i$ in $H$ are given by the $(q-1)$-line $L$, $\tfrac{q-1}{k}$ lines of multiplicity $0$, and 
	$\left(q-\tfrac{q-1}{k}\right)$ lines of multiplicity $k$. Now we are ready to determine the spectrum $(a_i)$ of $H$:
	\begin{eqnarray*}
		a_0 &=& 2\cdot\frac{q-1}{k} \\ 
		a_k &=& (q-1)\cdot (q-1)+2\cdot\left(q-\frac{q-1}{k}\right)\\ 
		a_{k+1} &=& (q-1)\cdot 1\\ 
		a_{q-1} &=& 1
	\end{eqnarray*} 
	and $a_i=0$ for all $i\notin\left\{0,k,k+1,q-1\right\}$. If $k=q-1$ or $k+1=q-1$, then we would have to take the sum of both 
	values, but we suppress this technical subtlety for the ease of notation. From the third standard equation 
	we conclude
	$$
	0={k\choose 2}a_k+{{k+1}\choose2} a_{k+1}+{{q-1}\choose 2}a_{q-1}-{{kq}\choose 2}=q\cdot \frac{(k - 1)(k - q + 1)}{2},
	$$
	so that $k\in\{1,q-1\}$. 
	%% (Note that the first and the second standard equation are satisfied, as it should be the case if the counting is done correctly.) 
	
	So, considering the $q+1$ hyperplanes through $L$ in $\PG(3,q)$ we conclude that $q$ have multiplicity $q$ and one has multiplicity $q(q-1)$, 
	where the latter contradicts \Cref{cor_cy_no_qm1_plane}. 
\end{Proof}

We remark that the preceding result can be used to simplify the proof of the cylinder conjecture for $q = 5$, although it is not necessary.

\begin{Lemma}
	\label{lemma_cy_no_qm3_plane}
	For a set of $q(q-3)$ points in $\PG(2,q)$, $q \geq 7$, there exists a line with a multiplicity not in $\{0,q-3,q-2\}$.
\end{Lemma}
\begin{Proof}
	Otherwise the first two standard equations would give  
	$a_{q-2}=a_0(q-3)-(q-3)$ and $a_{q-3} =q^2+2q-2 -a_0(q-2)$  
	so that the third yields
	$$
	0=a_{q-3}{{q-3}\choose 2}+a_{q-2}{{q-2}\choose 2}-{{q(q-3)}\choose 2}= %% wrong: \frac{(q-3)(a_0(q-2)+2q^2-2q+2)}{2},
	\frac{(q-3)(a_0(q-2)-4q+2)}{2}
	$$ 
	which implies
	$$
	a_0=\frac{4q-2}{q-2}=4+\frac{6}{q-2}.
	$$
	Since $a_0\in\N_0$, we have $q\in\{3,4,5,8\}$. Due to our assumption $q\ge 7$ it remains to exclude the case $q=8$, where $a_0=5$, $a_5=48$, and $a_6=20$. 
	Consider a $6$-line $L$. The only possibility for the distribution of the multiplicities of the lines through a $0$-point on $L$ is given by 
	$0^2 5^2 6^5$, so that there are $3\cdot 2=6$ $0$-lines. This contradicts $a_0=5$.
\end{Proof}
The assumption $q\ge 7$ in \Cref{lemma_cy_no_qm3_plane} is necessary, see \Cref{remark_10_0_2_3_arc_pg_2_5} for a counterexample for $q=5$ and the generator matrix
$$ 
\begin{pmatrix}
1&1&0&0\\
1&0&1&0\\
1&0&0&1
\end{pmatrix}  
$$
for $q=4$.

\begin{Corollary}
	\label{cor_cy_no_qm3_plane}
	Under $(\star)$, there can be no $q(q-3)$-plane.
\end{Corollary} 
\begin{Proof}
	Due to \Cref{cor_cy_conj_true_2_3}, \Cref{cor_cy_conj_true_4}, and \Cref{prop_cy_conj_true_5} we can assume $q\ge 7$. 
	\Cref{lemma_cy_gamma_2_at_most_q_minus_2} implies that the maximum line multiplicity is at most $q-2$. If $H$ is a $q(q-3)$-plane, then 
	\Cref{lemma_cy_line_multiplicities_in_plane} and \Cref{lemma_cy_no_qm3_plane} yield a contradiction.   
\end{Proof}

%\begin{Lemma}
%	\label{lemma_cy_no_qm4_plane}
%	For a set of $q(q-4)$ points in $\PG(2,q)$ with line multiplicities contained in $\{0,q-2,q-3,q-4\}$, we need $q\le 13$.
%\end{Lemma}
%\begin{Proof}
%	Solving the standard equations gives
%	\begin{eqnarray*}
%		a_{q-2} &=& \frac{(q-4)(5q-3)}{2}-\frac{(q-3)(q-4)}{2} a_0 \\ 
%		a_{q-3} &=& -(q-4)(5q-2)+(q-2)(q-4)a_0 \\ 
%		a_{q-4} &=& \frac{7q^2-19q+6}{2}  -\frac{(q-2)(q-3)}{2}a_0, 
%	\end{eqnarray*}
%	so that $a_{q-2}\ge 0$ implies 
%	$$
%	a_0\le \frac{5q-3}{q-3}= 6-\frac{q-15}{q-3}
%	$$
%	and 
%	$a_{q-3}\ge 0$ implies
%	$$
%	a_0 \ge \frac{5q-2}{q-2}=5+\frac{8}{q-2}.
%	$$
%	So, for $q\ge 16$ we have $a_0\notin \N_0$, which is a contradiction.
%\end{Proof}
%
%\begin{Remark}
%	\label{remark_cy_no_qm4_plane}
%	In the proof of \Cref{lemma_cy_no_qm4_plane} the only possible choices for $a_0$ are given by 
%	\begin{itemize}
%		\item $q=7$: $a_0\in\{7,8\}$, and $\left(a_0,a_{q-4},a_{q-3},a_{q-2}\right)=(7,38,6,6)$ or $(8,28,21,0)$;
%		\item $q=8$: $a_0=7$, and $\left(a_0,a_{q-4},a_{q-3},a_{q-2}\right)=(7,46,16,4)$; 
%		\item $q=9$: $a_0=7$, and $\left(a_0,a_{q-4},a_{q-3},a_{q-2}\right)=(7,54,30,0)$;
%		\item $q=11$: $a_0=6$, and $\left(a_0,a_{q-4},a_{q-3},a_{q-2}\right)=(6,106,7,14)$;
%		\item $q=13$: $a_0=6$, and $\left(a_0,a_{q-4},a_{q-3},a_{q-2}\right)=(6,141,27,9)$.
%	\end{itemize}  
%\end{Remark}

\begin{Lemma}
	\label{lemma_cy_no_qm4_plane}
	For a set of $q(q-4)$ points in $\PG(2,q)$, $q \geq 7$, with line multiplicities contained in $\{0,q-2,q-3,q-4\}$, we need $q = 7$ and spectrum $\left(a_0,a_{q-4},a_{q-3},a_{q-2}\right)=(8,28,21,0)$.
\end{Lemma}
\begin{Proof}
	Solving the standard equations gives
	\begin{eqnarray*}
		a_{q-2} &=& \frac{(q-4)(5q-3)}{2}-\frac{(q-3)(q-4)}{2} a_0 \\ 
		a_{q-3} &=& -(q-4)(5q-2)+(q-2)(q-4)a_0 \\ 
		a_{q-4} &=& \frac{7q^2-19q+6}{2}  -\frac{(q-2)(q-3)}{2}a_0, 
	\end{eqnarray*}
	so that $a_{q-2}\ge 0$ implies 
	$$
	a_0\le \frac{5q-3}{q-3}= 6-\frac{q-15}{q-3}
	$$
	and 
	$a_{q-3}\ge 0$ implies
	$$
	a_0 \ge \frac{5q-2}{q-2}=5+\frac{8}{q-2}.
	$$
	So, for $q\ge 16$ we have $a_0\notin \N_0$, which is a contradiction.
	
	For $7 \leq q \leq 13$ we can see that $a_{q-3} > 0$. So consider a $(q-3)$-line $L$. Through each of the four $0$-points on $L$ there go at least two $0$-lines, since otherwise $1\cdot 0+(q-1)\cdot (q-4)+1\cdot (q-3)=q(q-4)+1>|\cK|$. Therefore, we have $a_0 \geq 8$, which is only possible for $q = 7$. The spectrum then follows from the equations above.
\end{Proof}

%\begin{Lemma}
%	\label{lemma_cy_21_plane_q_7_spectrum}
%	A set $\cK$ of $21$ points in $\PG(2,7)$, whose line multiplicities are contained in $\{0,3,4,5\}$ has spectrum given by $a_0=8$, $a_3=28$, $a_4=21$, and $a_i=0$ otherwise.  
%\end{Lemma}
%\begin{Proof}
%	Given the proof of \Cref{lemma_cy_no_qm4_plane} it remains to exclude the possible spectrum $a_0=7$, $a_3=38$, $a_4=6$, 
%	and $a_5=6$. Consider a $4$-line $L$. Through each of the four $0$-points on $L$ there go at least two $0$-lines, since 
%	otherwise $1\cdot 0+6\cdot 3+1\cdot 4=22>21=|\cK|$. However, this gives at least $4\cdot 2=8>7=a_0$ $0$-lines, which is a 
%	contradiction.   
%\end{Proof}

%\begin{Remark}
%	A set of points with parameters and spectrum as specified in \Cref{lemma_cy_21_plane_q_7_spectrum} indeed exists.
%\end{Remark}

\begin{Lemma}
	\label{lemma_cy_q_7_no_21_plane}
	Under $(\star)$ for $q = 7$, the spectrum is given by $a_0=22$, $a_7=357$, $a_{14}=21$ and $a_{21} = 0$.
\end{Lemma}
\begin{Proof}
	First we will show that $a_{21} = 0$. Assume to the contrary the existence of a $21$-plane $H$. 	From \Cref{lemma_cy_no_qm4_plane} we conclude that the spectrum $\left(b_i\right)$ of $\cS \cap H$ satisfies $b_0=8$, $b_3=28$, $b_4=21$, and $b_i=0$ otherwise. Consider the possible hyperplane distributions through  
	a $0$-line in $H$: $0^5 7^1 21^2$, $0^5 14^2 21^1$, $0^4 7^2 14^1 21^1$, and $0^3 7^4 21^1$. In each case there are at least three $0$-planes through a $0$-line  
	in $H$, so that there are at least $3b_0=24$ $0$-planes in total. However, solving the standard equations for $\left\{a_0,a_7,a_{14}\right\}$ gives 
	$a_0 = 22 - a_{21}$, $a_7 = 357 + 3a_{21}$, and $a_{14} = 21 - 3a_{21}$, so that $a_0\le 22$, a contradiction. Therefore $a_{21} = 0$ and the values for $a_0, a_7$ and $a_{14}$ follow immediately as well.
\end{Proof}

\begin{Remark}
	\label{remark_gcc_no_line_mult_1_q_7}
	After 1671~seconds of computation time, \texttt{QextNewEdition} claims that there no sets of $14$ points with line multiplicity at most $5$ in $\PG(2,7)$ without 
	lines of multiplicity $1$. If we additionally assume that there is no line of multiplicity $5$, then 908~seconds of computation time are 
	needed. Using \Cref{lemma_cy_q_7_no_21_plane} this implies the truth of the cylinder conjecture for $(4,1,7)$.
\end{Remark} 

In the following we want to give an alternative, computer-free proof of the cylinder conjecture for $(4,1,7)$.

\begin{Lemma}
	\label{lemma_spectra_14_le5_arc_pg_2_7}
	A set $\cK$ of $14$ points in $\PG(2,7)$, whose line multiplicities are contained in $\{0,2,3,4,5\}$ has spectrum $\left(a_0,a_2,a_3,a_4,a_5\right)$ either $(10,37,2,8,0)$, $(11,31,10,5,0)$, $(12,25,18,2,0)$, 
	$(11,30,13,2,1)$, or $(10,36,5,5,1)$. 
\end{Lemma}
\begin{Proof}
	Solving the standard equations for $\left\{a_0,a_2,a_3\right\}$ gives
	\begin{eqnarray*}
		a_0 &=& \frac{38}{3} - a_5 - \frac{1}{3} a_4 \\ 
		a_2 &=& 21 + 5a_5 + 2a_4 \\ 
		a_3 &=& \frac{70}{3} - 5a_5 - \frac{8}{3} a_4 
	\end{eqnarray*} 
	From $a_0\in\N_0$ we conclude $a_4\equiv 2\pmod 3$, so that especially $a_4\ge 2$. With this, $a_3\in\N_0$ yields $a_4\le 8$ and $a_5\le 3$.
	
	In order to show $a_5\le 1$ we consider a $5$-line $L$. Now, let $P$ be an arbitrary $1$-point on $L$ and $Q$ be an arbitrary $0$-point on $L$. 
	Counting the points on the lines $L,L_1,\dots,L_7$ through $P$ gives
	$$
	14=|\cK|=\cK(L)+\sum_{i=1}^7 \cK(L_i)-7\cdot\cK(P)=\cK(L_1)+\sum_{i=2}^7 \cK(L_i)-2\ge \cK(L_1)+6\cdot 2 -2=\cK(L_1)+10,
	$$
	so that there is no $5$-line through $P$ besides $L$. Now assume that $M$ is a $5$-line through $Q$ and let $R$ be a $0$-point on $M$ 
	not equal to $Q$. By $L_0',\dots,L_7'$ we denote the lines through $R$, where we assume $L_0'=M$. Since five of the lines $L_1',\dots, L_7'$ 
	hit $L$ in a point they have multiplicity at least $2$, which yields at least $5+5\cdot 2=15>14$ points in $\cK$, which is a contradiction. 
	Thus, there is also no  $5$-line through $Q$ besides $L$, so that $a_5$ is at most $1$.  
	
	To sum up, if $a_5=0$, then $a_3\in\N_0$ implies $a_4\in\{2,5,8\}$, and if $a_5=1$, the  $a_3\in\N_0$ implies $a_4\in\{2,5\}$. Plugging this into the above 
	equations gives the five stated spectra.
\end{Proof}
Note that we have applied the same {\lq\lq}technique{\rq\rq} to conclude $a_5\le 1$ as the one used in the proof of 
\Cref{lemma_cy_gamma_2_at_most_q_minus_2}. We will now rule out all possibilities from \Cref{lemma_spectra_14_le5_arc_pg_2_7}.

%% \begin{Remark}
%%   The linear programming method applied to the stated equations gives only $a_5\le 7$. Restricting to integral solutions we obtain 
%%   $a_5\le 3$, which is attained for $a_0=33$, $a_2=13$, $a_3=0$, $a_4=8$, and $a_5=3$. Maximizing $a_4$ the LP method gives $a_4\le 13$ and 
%%   the ILP method gives $a_4\le 8$ attained at $a_0=24$, $a_2=16$, $a_3=9$, $a_4=8$, and $a_5=0$.
%% \end{Remark}

%\begin{Lemma}
%	\label{lemma_14_le5_arc_pg_2_7_no_8_0}
%	A set $\cK$ of $14$ points in $\PG(2,7)$, whose line multiplicities are contained in $\{0,2,3,4,5\}$ can not have spectrum $(10,37,2,8,0)$.
%\end{Lemma}
%\begin{Proof}
%	Assume that $\left(a_4,a_5\right)=(8,0)$, so that $a_0=10$, $a_2=37$, and $a_3=2$. 
%	Let $L$ be a $4$-line and $P$ be a $1$-point on $L$. Since there is no $5$-line at all, at least one $3$-line must go through $P$. (The possible 
%	distributions of the line multiplicities are $4^1 3^3 2^4$ and $4^2 3^1 2^5$.) This gives $a_3\ge 4$, which contradicts $a_3=2$.
%\end{Proof}

\begin{Lemma}
	\label{lemma_14_le5_arc_pg_2_7_a5_eq_0}
	A set $\cK$ of $14$ points in $\PG(2,7)$, whose line multiplicities are contained in $\{0,2,3,4,5\}$, can not have $5$-lines.   
\end{Lemma}
\begin{Proof}
	Let $L$ be a $5$-line and $Q_1$, $Q_2$, and $Q_3$ be the $0$-points on $L$. Define a new set $\cK'$ as the symmetric difference of $\cK$ and $L$: $\cK' = \cK \setminus (L \cap \cK) \cup \{Q_1,Q_2,Q_3\}$.
	As $\cK$ did not have any $1$-lines, we can see that $\cK'$ is a non-trivial blocking set of size $12$ in $\PG(2,7)$. 
	Blocking sets of cardinality $12$ in $\PG(2,7)$ not containing a line have been classified in \cite{blokhuis2003blocking}: besides the projective triangle there exists a unique sporadic example of non-R\'edei type. Both of these constructions share the property that there exists at least one $2$-line through every $1$-point. Now considering the point $Q_1$, we can easily see that there are no $2$-lines to $\cK'$ through it, as they would have to come from $1$-lines to $\cK$, which do not exist.
\end{Proof}

\begin{Lemma}
	\label{lemma_14_le5_arc_pg_2_7_no_5_0}
	No set $\cK$ of $14$ points in $\PG(2,7)$ with spectrum $\left(a_0,a_2,a_3,a_4\right)=(11,31,10,5)$ and $a_i = 0$ otherwise exists.
\end{Lemma}
\begin{Proof}
	The three possible distributions of the multiplicities of the lines through a $1$-point are
	$4^2 3^1 2^5$, $4^1 3^3 2^4$, $3^5 2^3$ and we speak of type $A_1$, $A_2$, and $A_3$, respectively.
	
	Assume that $P$ is a $1$-point of type $A_3$ and let $L$ be an arbitrary $3$-line through $P$. The five $0$-points
	on $L$ are contained in another $3$-line besides $L$, by parity considerations. Thus by $a_3 = 10$, the two $1$-points on $L$ that are not equal to $P$ are
	of type $A_1$. Since $L$ was chosen arbitrary, all ten $1$-points on $3$-lines through $P$ that are not equal to $P$
	are of type $A_1$. Let $Q_1$, $Q_2$, and $Q_3$ denote the three other $1$-points not equal to $P$ and not forming a
	$3$-line with $P$. All of the five $3$-lines not incident with $P$ have to consist of $1$-points in
	$\left\{Q_1,Q_2,Q_3\right\}$, which is impossible. Thus, there is no $1$-point of type $A_3$.
	
	Consider an arbitrary $3$-line $L'$. The five $0$-points on $L'$ are contained in another $3$-line besides $L'$, for the same reason as before.
	Thus not all three $1$-points on $L'$ can be of type $A_2$ as $a_3 = 10$. 
	%Since $a_4=5$ also not all three $1$-points on $L'$ can be	of type $A_1$. 
	So each $3$-line contains at least one $1$-point of type $A_1$, which is then contained in no other $3$-line.
	This means that there are at least ten $1$-points of type $A_1$. Counting the $4$-lines, this gives $a_4\ge (10\cdot 2+4\cdot 1)/4=6>5$, which is a
	contradiction.
	
\end{Proof}

\begin{Proposition}
	\label{prop_cy_conj_true_7}
	The generalized cylinder conjecture is true for the triple $(4,1,7)$.
\end{Proposition}
\begin{Proof}
	Assume the existence of a counterexample $\cS$. Due to \Cref{lemma_cy_q_7_no_21_plane}, there exists a hyperplane $H$ so that $\cS \cap H$ 
	is a set of $14$ point whose line multiplicities are contained in $\{0,2,3,4,5\}$. From \Cref{lemma_spectra_14_le5_arc_pg_2_7}, \Cref{lemma_14_le5_arc_pg_2_7_a5_eq_0} and \Cref{lemma_14_le5_arc_pg_2_7_no_5_0} the 
	spectrum of $\cS \cap H$ is given by $\left(b_0,b_2,b_3,b_4,b_5\right)=(12,25,18,2,0)$. Note that every $j$-line in $H$, where $j\ge 2$, is contained 
	in $(j-1)$ planes of multiplicity $14$ and $(9-j)$ planes of multiplicity $7$ in $\PG(3,7)$. 
	Thus $\cS$ contains at least $b_3+2b_4=18+2\cdot 2=22$ planes of multiplicity $14$, which contradicts $a_{14}=21$.
\end{Proof}

While our computer-free proof of the cylinder conjecture for $(4,1,7)$ is rather lengthy, most parts are more or less systematic and might be generalized to 
larger field sizes. A big obstacle is that we cannot prove the truth of the observation in \Cref{remark_gcc_no_line_mult_1_q_7} directly. Actually, we do 
not have a complete proof of this specific nonexistence result sets of $14$ points in $\PG(2,7)$ without $1$-lines and just sailed around the remaining 
open case in the proof of \Cref{prop_cy_conj_true_7}. Maybe other methods are more suitable for this kind of problems in $\PG(2,q)$. Of course, 
allowing computer enumerations drastically reduces the length of the argumentation. Starting from \Cref{lemma_cy_line_multiplicities_in_plane} and 
\Cref{lemma_cy_line_multiplicities_in_plane} we can computationally exclude many possibilities for the restriction $\cS \cap H$ for a hyperplane $H$. In other words, the possible multiplicities for the weights $\cS(H)$ for the hyperplanes can be restricted by 
enumeration results for $3$-dimensional codes over $\F_q$. By considering a subcode of the $4$-dimensional projective code corresponding to $\cS$ we obtain 
a $q^2$-divisible $\left[q^2-1,3\right]_q$-code $C$ with a restricted set of weights that might also be enumerated computationally. For our example we remark that 
there are $54$ non-isomorphic $[48,3,\{21,28,35,42\}]_7$-codes and $46$ non-isomorphic $[48,3,\{28,35,42\}]_7$-codes. Moreover, the information that $\cS$ does 
not contain a full affine line restricts the possible residual codes of codewords in $C$. By that criterion $6$ of of the $46$ non-isomorphic $[48,3,\{28,35,42\}]_7$-codes 
can be excluded. Similar restrictions can arise from the previously mentioned classification of projective $[kq,3,\{(k-1)q+1,\dots, k(q-1),kq\}]_q$-codes. E.g., as also 
theoretically proven, all projective $[21, 3, \{15, 16, 17, 18, 21\}] 7$-codes do not contain codewords of weight $15$ or $16$. So, in the residual code 
of a codeword of weight $28$ in $C$ the weights $15$ and $16$ cannot occur. This excludes $12$ further codes. For the remaining twenty-eight $[48,3,\{28,35,42\}]_7$-codes 
we can computationally check whether an extension to a projective $[49,4,\{28,35,42\}]_7$-code exists. To this end we can utilize and ILP formulation and an 
ILP solver. We remark that the tightest ILP instance needed $1\,238\,996$~branch\&bound nodes and 28.75~hours of computation time.  
At the very least, this approach gives a computational verification of \Cref{prop_cy_conj_true_7}. 

\vspace{1em}

Let us finish with some conclusions for the generalized cylinder conjecture for $q=8$. Assume that $\cS$ is an $8$-divisible spanning set of $64$ points in 
$\PG(3,q)$ that is not a cylinder. From \Cref{cor_cy_no_qm1_plane}, \Cref{cor_cy_no_qm2_plane}, \Cref{cor_cy_no_qm3_plane} and
\Cref{lemma_cy_no_qm4_plane}, we conclude that the hyperplane multiplicities 
with respect to $\cS$ are contained in $\{0,8,16,24\}$.

Solving the standard equations for the spectrum $\left(a_i\right)$ of $\cS$ gives
\begin{eqnarray*}
	a_0 &=& 29 - a_{24}\\
	a_8 &=& 528 + 3a_{24}\\ 
	a_{16} &=& 28 - 3a_{24}, 
\end{eqnarray*} 
so that $a_0\le 29$. Now assume that $H$ is a $24$-plane and consider the spectrum $\left(b_i\right)$ of $\cS \cap H$. Solving the standard equations 
for $\left\{b_3,b_5,b_6\right\}$ gives 
\begin{eqnarray*}
	b_3 &=& 97 - 5b_0 -\frac{b_4}{3} \\ 
	b_5 &=& -69 + 9b_0 - b_4 \\ 
	b_6 &=& 45 -5b_0 + \frac{b_4}{3},
\end{eqnarray*}
so that $b_5\ge 0$ implies $b_0\ge \left\lceil\tfrac{69}{9}\right\rceil=8$. Since through every $0$-line in $H$ there are at least three $0$-planes, 
$a_0\le 29$ implies $b_0\le \left\lfloor \tfrac{29}{3}\right\rfloor=9$. For $b_0=8$ we have $b_3 = 62 - b_6$, $b_4 = -15 + 3b_6$, and $b_5 = 18 - 3b_6$, 
so that either
$$
\left(b_0,b_3,b_4,b_5,b_6\right) \in\left\{(8,57,0,3,5),(8,56,3,0,6)\right\}
$$
or
$$
\left(b_0,b_3,b_4,b_5,b_6\right)=\left(9,52 - b_6,3b_6,12 - 3b_6,b_6\right),
$$
where $0\le b_6\le 4$. Consider a $4$-line $L$. Through each of the five $0$-points on $L$ there are at least two incident $0$-lines, so that $b_0\ge 5\cdot 2=10$. 
Thus, we conclude
$$
\left(b_0,b_3,b_4,b_5,b_6\right) \in\left\{(8,57,0,3,5),(9,52,0,12,0)\right\}.
$$
For the second case consider a $1$-point $P$. Since all lines through $P$ have to be $3$- or $5$-lines, we have $|\cS \cap H| \equiv 1\pmod 2$, which is a 
contradiction. For the first case we consider a $5$-line $L$ and observe that the unique possibility for the distribution of the multiplicities of the 
lines through a $1$-point on $L$ is given by $3^7 5^1 6^1$. Thus, besides $L$, there remain eight $0$-lines, twenty-two $3$-lines, and two $5$-lines 
for the four $0$-points on $L$. The only possibility for a $0$-point, using only $0$-, $3$-, and $5$-lines, is $0^3 3^3 5^3$ for the distribution of the 
multiplicities of the incident lines. This case cannot occur four times, so that we finally conclude $a_{24}=0$, which leaves the unique spectrum  
$\left(a_0,a_8,a_{16}\right)=(29,528,28)$ for $\cS$. The above considerations are elementary and easy, but a bit ad hoc. As for $q=7$, we are again in 
a situation where it seems that we are missing the right tools to tackle the problem in an elegant way. Of course, it is very likely that the cylinder 
conjecture is true for $q=8$.

%% \section*{Acknowledgment}

%% \bibliography{bibliography}

\begin{thebibliography}{10}

\bibitem{ball2008graph}
S.~Ball.
\newblock On the graph of a function in many variables over a finite field.
\newblock {\em Des. Codes Cryptogr.}, 47(1-3):159--164, 2008.

\bibitem{blokhuis2003blocking}
A.~Blokhuis, A.~E. Brouwer, and H.~A. Wilbrink.
\newblock Blocking sets in {${\rm PG}(2,p)$} for small {$p$}, and partial
  spreads in {${\rm PG}(3,7)$}.
\newblock {\em Adv. Geom.}, (suppl.):S245--S253, 2003.
\newblock Special issue dedicated to Adriano Barlotti.

\bibitem{de2019cylinder}
J.~De~Beule, J.~Demeyer, S.~Mattheus, and P.~Sziklai.
\newblock On the cylinder conjecture.
\newblock {\em Des. Codes Cryptogr.}, 87(4):879--893, 2019.

\bibitem{govaerts2003particular}
P.~Govaerts and L.~Storme.
\newblock On a particular class of minihypers and its applications. {I}. {T}he
  result for general {$q$}.
\newblock {\em Des. Codes Cryptogr.}, 28(1):51--63, 2003.

\bibitem{partial_spreads_and_vector_space_partitions}
T.~Honold, M.~Kiermaier, and S.~Kurz.
\newblock Partial spreads and vector space partitions.
\newblock In {\em Network Coding and Subspace Designs}, pages 131--170.
  Springer, 2018.

\bibitem{kurz2019lincode}
S.~Kurz.
\newblock Lincode -- computer classification of linear codes.
\newblock {\em arXiv preprint 1912.09357}, 2019.

\bibitem{lavrauw2015fieldreduction}
M.~Lavrauw and G.~Van~de Voorde.
\newblock Field reduction and linear sets in finite geometry.
\newblock In {\em Topics in finite fields}, volume 632 of {\em Contemp. Math.},
  pages 271--293. Amer. Math. Soc., Providence, RI, 2015.

\bibitem{lovasz1983}
L.~Lov\'{a}sz and A.~Schrijver.
\newblock Remarks on a theorem of {R}\'{e}dei.
\newblock {\em Studia Sci. Math. Hungar.}, 16(3-4):449--454, 1983.

\bibitem{redei1970}
L.~R\'{e}dei.
\newblock {\em L\"{u}ckenhafte {P}olynome \"{u}ber endlichen {K}\"{o}rpern}.
\newblock Birkh\"{a}user Verlag, Basel-Stuttgart, 1970.

\bibitem{vandevoorde2011}
G.~Van~de Voorde.
\newblock On sets without tangents and exterior sets of a conic.
\newblock {\em Discrete Math.}, 311(20):2253--2258, 2011.

\bibitem{ward2001divisible}
H.~N. Ward.
\newblock Divisible codes---a survey.
\newblock {\em Serdica Math. J.}, 27(4):263--278, 2001.

\end{thebibliography}
%% %\bibdata{cylinder_conjecture}
%% %%\bibliographystyle{amsplain}
%% \bibliographystyle{abbrv}

\end{document}